\documentclass[a4paper,11pt]{amsart}

\usepackage{graphicx}
\usepackage{mathptmx}
\usepackage{amsmath}
\usepackage{amssymb}
\usepackage{enumitem}
\usepackage{xcolor}
\usepackage{pgfplots}

\newmuskip\pFqmuskip

\newcommand*\pFq[6][8]{%
  \begingroup 
  \pFqmuskip=#1mu\relax
  \mathcode`=\string"8000
  \begingroup\lccode`\~=`\,
  \lowercase{\endgroup\let~}\pFqcomma
  F^{#2}_{#3}{\left(\genfrac..{0pt}{}{#4}{#5}\bigg|#6\right)}%
  \endgroup
}
\newcommand{\pFqcomma}{\mskip\pFqmuskip}

\newtheorem{theorem}{Theorem}[section]

\begin{document}

\title[A note on new type degenerate Srirling numbers of the first kind]{A note on new type degenerate Srirling numbers of the first kind}

\author{Taekyun  Kim}
\address{Department of Mathematics, Kwangwoon University, Seoul 139-701, Republic of Korea}
\email{tkkim@kw.ac.kr}
\author{Dae San  Kim}
\address{Department of Mathematics, Sogang University, Seoul 121-742, Republic of Korea}
\email{dskim@sogang.ac.kr}

\author{Kyo-Shin Hwang}
\address{Graduate School of Education, Yeungnam University, Gyeongsan 38541, Republic of Korea}
\email{kshwang@yu.ac.kr}

\author{Dmitry V. Dolgy}
\address{Kwangwoon Global Education Center, Kwangwoon University, Seoul, Republic of Korea}
\email{d\_dol@kw.ac.kr}

\subjclass[2010]{11B73; 11B83; 60-08}
\keywords{new type degenerate Stirling numbers of the first kind; unsigned new type degenerate Stirling numbers of the first kind; new type degenerate Stirling numbers of the second kind}

\begin{abstract}
We introduce a new sequence of unsigned degenerate Stirling numbers of the first kind. Following the work of Adell-Lekuona, who represented unsigned Stirling numbers of the first kind as multiples of the expectations of specific random variables, we express our new numbers as finite sums of multiples of the expectations of certain random variables. We also provide a representation of these new numbers as finite sums involving the classical unsigned Stirling numbers of the first kind. As an inversion formula, we define a corresponding sequence of new type degenerate Stirling numbers of the second kind. We derive expressions for these numbers as finite sums that involve the Stirling numbers of the second kind.
\end{abstract}

\maketitle

\section{Introduction}
For any nonzero $\lambda \in \mathbb{R}$, the degenerate exponentials are defined by
\begin{equation}\label{1}
e_\lambda^x (t) = \sum_{k=0}^\infty (x)_{k,\lambda}\frac{t^k}{k!},\ \ e_\lambda (t) =e_\lambda^1 (t), \quad(\rm see\ [2,4-7,12]),
\end{equation}
where
\begin{displaymath}
(x)_{0,\lambda}=1,\quad (x)_{n,\lambda}=x(x-\lambda)(x-2\lambda)\cdots(x-(n-1)\lambda),\quad (n\ge 1).
\end{displaymath}
Note that
\begin{displaymath}
\lim_{\lambda\to 0} e_\lambda(t) =e^t.
\end{displaymath} \par
It is well known that Stirling  numbers of the first kind are defined by
\begin{equation}\label{2}
\frac{1}{k!}\log^k(1+t) = \sum_{n=k}^\infty S_1 (n, k) \frac{t^n}{n!},\quad(k\ge 0), \quad(\rm see \ [3,9]).
\end{equation}
The unsigned Stirling numbers of the first kind are defined by
\begin{displaymath}
{n \brack k}= (-1)^{n-k} S_1(n,k), \quad(n \ge k\ge 0).
\end{displaymath}
Thus, we note from \eqref{2}  that
 \begin{equation}\label{3}
\frac{1}{k!}\log^k \left( \frac1{1-t}\right) =\sum_{n=k}^\infty {n \brack k}
\frac{t^n}{n!}, \quad(\rm see\ [3,6]).
\end{equation}
As the inversion formula of \eqref{2}, the Stirling numbers of the second kind are defined by
\begin{equation}\label{4}
\frac1{k!} \left(e^t-1\right)^k =\sum_{n=k}^\infty {n \brace k} \frac{t^n}{n!},\quad(n \ge k\ge 0), \quad(\rm see\ [3,9,11]).
\end{equation}

Let $\log_\lambda (t)$ be the degenerate logarithm, which is the compositional inverse of $e_\lambda(t)$. Then we note that
$$\log_\lambda(1+t)=\sum_{n=1}^\infty \lambda^{n-1}(1)_{n,\frac{1}{\lambda}}\frac{t^n}{n!} = \frac1\lambda\big((1+t)^\lambda -1\big),\quad(\rm see\ [6, 12]).$$
Recently, the degenerate Stirling numbers of the first kind are given by
\begin{equation*}
\frac1{k!} \log_\lambda^k (1+t) = \sum_{n=k}^\infty S_{1,\lambda}(n,k) \frac{t^n}{n!}, \quad(\rm see\ [6,12]).
\end{equation*}
The unsigned degenerate Stirling numbers of the first kind are defined by
\begin{displaymath}
{n \brack k}_{\lambda}= (-1)^{n-k}S_{1,\lambda}(n,k), \quad (n,k\ge 0).
\end{displaymath}
Thus we have
\begin{equation*}
\frac{1}{k!} \left( -\log_\lambda(1-t)\right)^k =\frac{1}{k!}\log_{-\lambda}^{k}\Big(\frac{1}{1-t}\Big)= \sum_{n=k}^\infty {n \brack k}_{\lambda}\frac{t^n}{n!}, \quad(\rm see\ [6,12]).
\end{equation*}
The degenerate Stirling numbers of the second kind are defined by
\begin{equation*}
\frac{1}{k!} \left( e_\lambda (t)-1\right)^k = \sum_{n=k}^\infty S_{2,\lambda}(n,k) \frac{t^n}{n!}, \quad(\rm see\ [6,12]).
\end{equation*} \par
Let $U$ and $X$ be two independent random variables where $U$ is the uniform random variable on $(0,1)$ and
$X$ is the exponential random variable with parameter $1$. Recall that the probability density function of $X$ is given by (see [10])
\begin{displaymath}
f_X (x) = \begin{cases}e^{-x}, &\text{if}\ \ x\ge0,\\ 0, &\text{if}\ \ x<0, \end{cases}
\end{displaymath}
and that the probability density function of $U$ is given by
\begin{displaymath}
g_U (x) = \begin{cases} 1, &\text{if}\ \ x\in (0,1),\\ 0, &\text{if}\ \ x\notin (0,1). \end{cases}
\end{displaymath}
Let $(U_j)_{j\ge 1}$ and $(X_j)_{j\ge1}$ be two sequences of independent copies of $U$ and $X$, respectively, both of them
mutually independent. We use the notation
\begin{displaymath}
S_k = U_1 X_1 + U_2X_2 + \cdots+U_kX_k, \ k=1,2, \cdots , \quad S_0 =0.
\end{displaymath} \par
Adell-Lekuona [1] showed the following identity (see \eqref{2}, \eqref{3}):
\begin{equation}\label{5}
{n \brack k}=\binom{n}{k}E\left[S_k^{n-k}\right], \quad (n \ge k \ge 1),
\end{equation}
where $E$ is the mathematical expectation. \\
Their idea of proof is to note that $E\left[e^{t S_{k}}\right]=\Big(\frac{\log(1-t)}{-t}\Big)^{k}$. \\
The aim of this paper is to derive a degenerate version of \eqref{5}. Namely, we show the following expression in Theorem 2.2:
\begin{equation}\label{6}
{n \brack k}_{-\lambda}^{*} =\sum_{m=k}^n \lambda^{m-k} \binom{n}{m} S_1(m,k) E\left[ (S_k)_{n-m, \lambda}\right],\quad (n \ge k \ge 1),
\end{equation}
where ${n \brack k}_{-\lambda}^{*}$ are the unsigned new type degenerate Stirling numbers of the first kind given by $\frac1{k!} \log^k \left(\frac1{1-\frac1\lambda \log(1+\lambda t)}\right)=\sum_{n=k}^\infty {n \brack k}_{-\lambda}^{*}\frac{t^n}{n!}$,\,\, (see \eqref{1}, \eqref{8}). \\
Here our idea is to note $E\left[e_{\lambda}^{S_{k}}(t)\right]= \Big(\frac1{\frac1\lambda \log(1+\lambda t)}\Big)^k \log^k \Big(\frac1{1- \frac1\lambda \log(1+\lambda t)}\Big)$. In this way, we were led to introduce the unsigned new type degenerate Stirling numbers of the first kind. Observe that \eqref{6} boils down to \eqref{5} if we let $\lambda \rightarrow 0$. In Theorem 2.1, we show that ${n \brack k}_{\lambda}^{*} = \sum_{m=k}^n \lambda^{n-m}{n \brack m}{m \brack k},\,\, (n \ge k \ge 0)$. Using this we compute ${n \brack k}_{\lambda}^{*}$, for $0 \le n \le 6$. As an inversion formula, we define a corresponding sequence ${n \brace k}_{\lambda}^{*}$, called the new type degenerate Stirling numbers of the second kind (see \eqref{17}). We show in Theorem 3.1 that ${n \brace k}_{\lambda}^{*}=\sum_{m=k}^n \lambda^{m-k} {m \brace k} {n \brace m},\,\,(n \ge k \ge 0)$. Using this, we compute ${n \brace k}_{\lambda}^{*}$, for $0 \le n \le 6$. As general references for this paper, the reader may refer to [3,9,10]. \par

\section{A new type degenerate Stirling numbers of the first kind}
For any nonzero $\lambda \in \mathbb{R}$, we consider the {\it{new type degenerate Stirling numbers of the first kind}} defined by
\begin{equation}\label{7}
\frac{1}{k!} \log^k  \Big(1+ \frac1\lambda \log(1+\lambda t)\Big)
= \sum_{n=k}^\infty S_{1,\lambda}^{*}(n,k) \frac{t^n}{n!},\quad(k\ge 0).
\end{equation}
Note that (see \eqref{2})
\begin{displaymath}
\lim_{\lambda\to 0} S_{1,\lambda}^{*} (n,k) = S_{1}(n,k).
\end{displaymath}
In addition, we define the {\it{unsigned new type degenerate Stirling numbers of the first kind}} by
\begin{displaymath}
{n \brack k}_{\lambda}^{*} =(-1)^{n-k} S_{1,\lambda}^{*} (n,k),\quad (n,k\ge 0).
\end{displaymath}
Then we note from \eqref{7} that
\begin{equation}\label{8}
\frac{1}{k!} \log^k  \bigg( \frac1{1+\frac1{\lambda} \log(1-\lambda t)}\bigg)
 = \sum_{n=k}^\infty {n \brack k}_{\lambda}^{*} \frac{t^n}{n!}.
\end{equation}
Note that (see \eqref{3})
\begin{equation} \label{9}
\lim_{\lambda\to 0} {n \brack k}_{\lambda}^{*}= {n \brack k},\quad(n \ge k\ge 0).
\end{equation} \par
From \eqref{8}, we note that
\begin{equation}\label{10}
\begin{split}
\sum_{n=k}^\infty {n \brack k}_{\lambda}^{*} \frac{t^n}{n!}
&= \frac1{k!} \log^k \bigg(\frac1{1+\frac1\lambda \log(1-\lambda t)}\bigg)\\
&= \sum_{m=k}^\infty {m \brack k}(-1)^{m} \lambda^{-m}\frac1{m!} \log^m (1-\lambda t)\\
&=\sum_{m=k}^\infty {m \brack k} (-1)^{m}\lambda^{-m}\sum_{n=m}^\infty S_1 (n,m) (-1)^{n}\lambda^n \frac{t^n}{n!}\\
&=\sum_{n=k}^\infty \sum_{m=k}^n \lambda^{n-m}{n \brack m}{m \brack k} \frac{t^n}{n!}.
\end{split}
\end{equation}
Therefore, by \eqref{10}, we obtain the following theorem.\par

\begin{theorem}
For any integers $n,k$ with $n \ge k\ge 0$, we have
\begin{equation*}
{n \brack k}_{\lambda}^{*} = \sum_{m=k}^n \lambda^{n-m}{n \brack m}{m \brack k}.
\end{equation*}
\end{theorem}

Using Theorem 2.1, we illustrate the values of the new type degenerate Stirling numbers of the first kind ${n \brack k}_{\lambda}^{*}$, for $n \le 6$ in the following.
We observe first that ${n \brack n}_{\lambda}^{*}$=1, for any nonnegative integer $n$; ${n \brack 0}_{\lambda}^{*}=0$, for any positive integer $n$, and ${0 \brack 0}=1$; ${n \brack n-1}_{\lambda}^{*}={n \brack n-1}(\lambda+1)=\binom{n}{2}(\lambda +1)$. We note that, as a plynomial in $\lambda$, the leading coefficient and the constant term of ${n \brack k}_{\lambda}^{*} = \sum_{m=k}^n \lambda^{n-m}{n \brack m}{m \brack k}$ are the same ${n \brack k}$.
\vspace{0.1in}
\begin{align*}
&{2 \brack 1}_{\lambda}^{*}=\lambda+1,
{3 \brack 1}_{\lambda}^{*}=2 \lambda^{2}+3 \lambda+2,
{3 \brack 2}_{\lambda}^{*}=3 \lambda +3, \\
&{4 \brack 1}_{\lambda}^{*}=6 \lambda^{3}+11 \lambda^{2}+12 \lambda +6,
{4 \brack 2}_{\lambda}^{*}=11 \lambda^{2}+18 \lambda +11,
{4 \brack 3}_{\lambda}^{*}=6 \lambda +6, \\
&{5 \brack 1}_{\lambda}^{*}=24 \lambda^{4}+50 \lambda^{3}+70 \lambda^{2}+60 \lambda +24,
{5 \brack 2}_{\lambda}^{*}=50 \lambda^{3}+105 \lambda^{2}+110 \lambda +50, \\
&{5 \brack 3}_{\lambda}^{*}=35 \lambda^{2}+60 \lambda +35,
{5 \brack 4}_{\lambda}^{*}=10 \lambda +10, \\
&{6 \brack 1}_{\lambda}^{*}=120 \lambda^{5}+274 \lambda^{4}+450 \lambda^{3}+510 \lambda^{2}+360 \lambda +120, \\
&{6 \brack 2}_{\lambda}^{*}=274 \lambda^{4}+675 \lambda^{3}+935 \lambda^{2}+750 \lambda +274, \\
&{6 \brack 3}_{\lambda}^{*}=225 \lambda^{3}+510 \lambda^{2}+525 \lambda+225,
{6 \brack 4}_{\lambda}^{*}=85 \lambda^{2}+150 \lambda+85,
{6 \brack 5}_{\lambda}^{*}=15 \lambda+15.
\end{align*} \par
\vspace{0.1in}
Now, we observe that
\begin{equation}\label{11}
\begin{split}
E\left[ e_\lambda^{UX}{(t)}\right] &= \int_0^1\int_0^\infty e^{-x(1-\frac{u}{\lambda}\log(1+\lambda t))} dx du\\
&=\int_0^1 \frac1{1-\frac{u}{\lambda}\log(1+\lambda t)}du\\
&= - \frac{\log\left(1-\frac{1}{\lambda}\log(1+\lambda t)\right)}{\frac{1}{\lambda}\log(1+\lambda t)}
= \frac{\log\left(\frac1{1-\frac{1}{\lambda}\log(1+\lambda t)}\right)}{\frac{1}{\lambda}\log(1+\lambda t)}.
\end{split}
\end{equation}
From \eqref{11}, we note that
\begin{equation}\label{12}
\begin{split}
E\left[ e_\lambda^{S_k }{(t)}\right] &=  E\left[ e_\lambda^{U_1 X_1 }{(t)}\right] E\left[ e_\lambda^{U_2 X_2 }{(t)}\right]\cdots
 E\left[ e_\lambda^{U_k X_k }{(t)}\right]\\
&= \bigg(\frac1{\frac1\lambda \log(1+\lambda t)}\bigg)^k \log^k \bigg(\frac1{1- \frac1\lambda \log(1+\lambda t)}\bigg).
\end{split}
\end{equation}
By \eqref{12}, we get
\begin{equation}\label{13}
\begin{split}
E\left[ e_\lambda^{S_k }{(t)}\right]\frac1{k!} \left(\frac1\lambda \log(1+\lambda t)\right)^k
&= \frac1{k!} \log^k \left(\frac1{1-\frac1\lambda \log(1+\lambda t)}\right)\\
&=\sum_{n=k}^\infty {n \brack k}_{-\lambda}^{*}\frac{t^n}{n!}.
\end{split}
\end{equation} \par
On the other hand, by \eqref{1} and \eqref{2}, we get
\begin{equation}\label{14}
\begin{split}
&E\left[ e_\lambda^{S_k }{(t)}\right] \frac1{k!} \left(\frac1\lambda \log(1+\lambda t)\right)^k\\
&=\sum_{l=0}^\infty E\left[ (S_k)_{l, \lambda}\right]  \frac{t^l}{l!} \sum_{m=k}^\infty \lambda^{m-k} S_1(m,k)\frac{t^m}{m!} \\
&=\sum_{n=k}^\infty \sum_{m=k}^n\binom{n}{m}  E\left[ (S_k)_{n-m, \lambda}\right] \lambda^{m-k} S_1(m,k)\frac{t^n}{n!}.
\end{split}
\end{equation}
Therefore, by \eqref{13} and \eqref{14}, we obtain the following theorem.\par

\begin{theorem}
For any integers $n,k$ with $n\ge k\ge 1$, we have
\begin{equation}\label{15}
{n \brack k}_{-\lambda}^{*} =\sum_{m=k}^n \lambda^{m-k} \binom{n}{m} S_1(m,k) E\left[ (S_k)_{n-m, \lambda}\right] .
\end{equation}
\end{theorem}
By taking $\lambda \rightarrow 0$ in \eqref{15} and using \eqref{9}, we recover the equation \eqref{5}
\begin{displaymath}
{n \brack k}= \lim_{\lambda\to 0} {n \brack k}_{-\lambda}^{*}
= \binom{n}{k}  E\left[ S_k^{n-k}\right], \quad(\rm see\ [1]).
\end{displaymath}
From \eqref{15}, we note that
\begin{equation}\label{16}
\begin{split}
{n \brack k}_{-\lambda}^{*}
&=\sum_{m=k}^n \binom{n}{m}  E\left[ (S_k)_{n-m, \lambda}\right]   \lambda^{m-k} S_1(m,k)\\
&= \sum_{m=k}^n \binom{n}{m}\lambda^{m-k} S_1 (m,k ) \underbrace{\int_0^1 \cdots\int_0^1}_{k-\text{times}}\\
&\qquad\times
\underbrace{\int_0^\infty\cdots \int_0^\infty}_{k-\text{times}} \Big(
\sum_{i=1}^k u_i x_i \Big)_{n-m,\lambda}e^{-(x_1+\cdots+x_k)}dx_1\cdots dx_k du_1 \cdots du_k.
\end{split}
\end{equation}
Therefore, by \eqref{16}, we obtain the following theorem.
\par
\begin{theorem}
For any integers $n,k$ with $n\ge k\ge 1$, we have
\begin{equation*}
\begin{split}
{n \brack k}_{-\lambda}^{*} &= \sum_{m=k}^n \binom{n}{m}\lambda^{m-k} S_1 (m,k ) \underbrace{\int_0^1 \cdots\int_0^1}_{k-\text{times}}\\
&\qquad\times
\underbrace{\int_0^\infty\cdots \int_0^\infty}_{k-\text{times}} \Big(
\sum_{i=1}^k u_i x_i \Big)_{n-m,\lambda}e^{-(x_1+\cdots+x_k)}dx_1\cdots dx_k du_1 \cdots du_k.
\end{split}
\end{equation*}
 \end{theorem}
\par

\section{Further Remark}
As the inversion formula of \eqref{7}, we define the {\it{new type degenerate Stirling numbers of the second kind}} by
\begin{equation}\label{17}
\frac{1}{k!} \left( \frac1\lambda \left(e^{\lambda (e^t -1)}-1\right)\right)^k
 = \sum_{n=k}^\infty {n \brace k}_{\lambda}^{*} \frac{t^n}{n!},\quad(k\ge 0).
\end{equation}
Note that (see \eqref{4})
\begin{displaymath}
\lim_{\lambda\to 0} {n \brace k}_{\lambda}^{*}= {n \brace k}, \quad(n \ge k\ge 0).
\end{displaymath}
From \eqref{17}, we have
\begin{equation}\label{18}
\begin{split}
\sum_{n=k}^\infty {n \brace k}_{\lambda}^{*} \frac{t^n}{n!}
&= \frac{1}{k!} \left( \frac1\lambda \left(e^{\lambda (e^t -1)}-1\right)\right)^k \\
&= \lambda^{-k}  \sum_{m=k}^\infty {m \brace k} \lambda^m  \frac1{m!} \left(e^t -1\right)^m\\
&= \lambda^{-k}  \sum_{m=k}^{\infty} {m \brace k} \lambda^m \sum_{n=m}^\infty {n \brace m}\frac{t^n}{n!} \\
& = \sum_{n=k}^\infty \sum_{m=k}^n \lambda^{m-k}{m \brace k}
{n \brace m} \frac{t^n}{n!}.
\end{split}
\end{equation}
Thus, by \eqref{18}, we obtain the following theorem.
\par

\begin{theorem}
For any integers $n,k$ with $n \ge k\ge 0$, we have
\begin{displaymath}
{n \brace k}_{\lambda}^{*}=\sum_{m=k}^n \lambda^{m-k} {m \brace k} {n \brace m}.
\end{displaymath}
\end{theorem} \par

Using Theorem 3.1, we illustrate the values of the unsigned new type degenerate Stirling numbers of the first kind ${n \brace k}_{\lambda}^{*}$, for $n \le 6$ in the following.
We observe first that ${n \brace n}_{\lambda}^{*}$=1, for any nonnegative integer $n$; ${n \brace 0}_{\lambda}^{*}=0$, for any positive integer $n$, and ${0 \brace 0}=1$; ${n \brace n-1}_{\lambda}^{*}={n \brace n-1}(\lambda +1)=\binom{n}{2}(\lambda +1)$, for any integer $n \ge 2$. Also, we note that, as a polynomial in $\lambda$, the leading coefficient and the constant term of ${n \brace k}_{\lambda}^{*}=\sum_{m=k}^n \lambda^{m-k} {m \brace k} {n \brace m}$ are the same number ${n \brace k}$.
\vspace{0.1in}
\begin{align*}
&{2 \brace 1}_{\lambda}^{*}=\lambda+1,
{3 \brace 1}_{\lambda}^{*}=\lambda^{2}+3 \lambda+1,
{3 \brace 2}_{\lambda}^{*}=3 \lambda +3, \\
&{4 \brace 1}_{\lambda}^{*}= \lambda^{3}+6 \lambda^{2}+7 \lambda +1,
{4 \brace 2}_{\lambda}^{*}=7 \lambda^{2}+18 \lambda +7,
{4 \brace 3}_{\lambda}^{*}=6 \lambda +6, \\
&{5 \brace 1}_{\lambda}^{*}=\lambda^{4}+10 \lambda^{3}+25 \lambda^{2}+15 \lambda +1,
{5 \brace 2}_{\lambda}^{*}=15 \lambda^{3}+70 \lambda^{2}+75 \lambda +15, \\
&{5 \brace 3}_{\lambda}^{*}=25 \lambda^{2}+60 \lambda +25,
{5 \brace 4}_{\lambda}^{*}=10 \lambda +10, \\
&{6 \brace 1}_{\lambda}^{*}= \lambda^{5}+15 \lambda^{4}+65 \lambda^{3}+90 \lambda^{2}+31 \lambda +1, \\
&{6 \brace 2}_{\lambda}^{*}=31\lambda^{4}+225 \lambda^{3}+455 \lambda^{2}+270 \lambda +31, \\
&{6 \brace 3}_{\lambda}^{*}=90 \lambda^{3}+375 \lambda^{2}+390 \lambda+90,
{6 \brace 4}_{\lambda}^{*}=65 \lambda^{2}+150 \lambda+65,
{6 \brace 5}_{\lambda}^{*}=15 \lambda+15.
\end{align*}
\vspace{0.1in}

\section{Conclusion}
From the identity $E\left[e^{t S_{k}}\right]=\Big(\frac{\log(1-t)}{-t}\Big)^{k}$, Adell-Lekuona derived the following identity:
\begin{equation*}
{n \brack k}=\binom{n}{k}E\left[S_k^{n-k}\right], \quad (n \ge k \ge 1).
\end{equation*}
In this paper, by using the identity $E\left[e_{\lambda}^{S_{k}}(t)\right]= \Big(\frac1{\frac1\lambda \log(1+\lambda t)}\Big)^k \log^k \Big(\frac1{1- \frac1\lambda \log(1+\lambda t)}\Big)$, we were able to deduce a degenerate version of Adell-Lekuona identity. Namely, we obtained
\begin{equation*}
{n \brack k}_{-\lambda}^{*} =\sum_{m=k}^n \lambda^{m-k} \binom{n}{m} S_1(m,k) E\left[ (S_k)_{n-m, \lambda}\right],\quad (n \ge k \ge 1).
\end{equation*}
This led us to the introduction of the unsigned new type degenerate Stirling numbers of the first kind. Furthemore, the explicit expression ${n \brack k}_{\lambda}^{*} = \sum_{m=k}^n \lambda^{n-m}{n \brack m}{m \brack k}$ was found for $n \ge k \ge 0$. As an inversion formula, we also defined a corresponding sequence ${n \brace k}_{\lambda}^{*}$ of new type degenerate Stirling numbers of the second kind. Then we show that ${n \brace k}_{\lambda}^{*}=\sum_{m=k}^n \lambda^{m-k} {m \brace k} {n \brace m},\,\,(n \ge k \ge 0)$. \par
In recent years, we have worked on probabilistic extensions of many special numbers and polynomials (see [7,8,12] and the references therein), and on applications of probabilty theory to the study of such numbers and polynomials (see [4] and the references therein). It is one of our research projects to continue to explore this line of research.
\par

\end{document}